\documentclass[reqno]{amsart}
\usepackage{amssymb}
\usepackage{graphics}
\usepackage{amscd}

\newtheorem{thm}{Theorem}
\newtheorem{cor}[thm]{Corollary}
\newtheorem{prop}[thm]{Proposition}
\newtheorem{lem}[thm]{Lemma}
\theoremstyle{definition}
\newtheorem{defn}{Definition}
\newtheorem{rem}[thm]{Remark}

\def\Alg{\mathrm{Alg}}
\def\Dialg{\mathrm{Dialg}}
\def\Var{\mathrm{Var}}
\def\idd{\mathop{\fam 0 id}\nolimits}

\def\Sym{\mathrm{Sym}}
\def\Span{\mathop{\fam 0 Span}\nolimits}
\def\Coeff{\mathop { \fam 0 Coeff}\nolimits}
\def\Ker{\mathop{\fam 0 Ker}\nolimits}
\def\Ress{\mathop {\fam 0 Res}}
\def\Cend{\mathop{\fam 0 Cend}\nolimits}
\def\End{\mathop{\fam 0 End}\nolimits}
\def\gc{\mathop{\fam 0 gc}\nolimits}
\def\Curr{\mathop{\fam 0 Cur}\nolimits}

\def\lbar{\dashv }
\def\rbar{\vdash }

\def\Comp{\mathrm{Comp}}
\def\Hom{\mathrm{Hom}}
\def\Vec{\mathrm{Vec}}
\def\Hmod{H\mbox{-}\mathrm{mod}}
\def\oo#1{\mathrel{{}_{(#1)}}}

\begin{document}

\title{Varieties of dialgebras and conformal algebras}

\thanks{Partially supported by RFBR 05--01--00230,
Complex Integration Program SB RAS (2006--1.9).
The author gratefully acknowledges the support from Deligne 2004 Balzan
prize in mathematics.}

\author{Pavel Kolesnikov}
\email{pavelsk@math.nsc.ru}

\address{Sobolev Institute of Mathematics,
Novosibirsk, Russia}

\begin{abstract}
For a given variety $\Var $ of algebras we define
the variety $\Var $ of dialgebras. This construction
turns to be closely related with varieties of pseudo-algebras:
every $\Var $-dialgebra can be embedded into an
appropriate pseudo-algebra of the variety $\Var $.
In particular, Leibniz algebras are exactly Lie dialgebras,
and every Leibniz algebra can be embedded into current
Lie conformal algebra.
\end{abstract}



\maketitle

\section*{Introduction}

Leibniz algebras, introduced by J.-L.~Loday \cite{L0}, are non-commutative
analogues of Lie algebras. The defining identity of the variety
of left (right) Leibniz algebras is the condition stating that
each operator of left (right) multiplication is a derivation, i.e.,
\[
x(yz)=(xy)z + y(xz)\quad \mbox{or}\quad
(xy)z = (xz)y + x(yz),
\]
respectively.
As Lie algebras are embeddable into associative algebras, Leibniz algebras
can be embedded into associative dialgebras.
By the definition \cite{L1, L2}, dialgebras are linear spaces with
two bilinear products $\lbar $, $\rbar $. The identities required
for a dialgebra $A$ to be associative are chosen in such a way that
the new operation
\[
ab=  a\rbar b - b\lbar a \quad \mbox{or}\quad  ab = a\lbar b - b\rbar a
\]
turns $A$ into a left (right) Leibniz algebra:
\begin{gather}
(x\lbar y)\rbar z = (x\rbar y)\rbar z,\quad
  x\lbar (y\rbar z)=x\lbar (y\lbar z),        \label{assD1} \\
(x\rbar y)\rbar z = x\rbar (y\rbar z),
   \quad  (x\lbar y)\lbar z = x\lbar (y\lbar z),  \label{assD2} \\
 (x\rbar y)\lbar z = x\rbar (y\lbar z).           \label{assD3}
\end{gather}

In the paper \cite{DL} the notion of an alternative dialgebra was
introduced.
 It was noted in
\cite{Faulk} that an algebra $A$ is embedded into the Steinberg Lie
algebra $\mathrm {st}(n,A)$, $n\ge 4$, if and only if $A$ is associative;
and if $n=3$ then $A$ has to be at least alternative. In \cite{DL},
the Steinberg Leibniz algebra $\mathrm{stl}(n,D)$ was introduced,
where $D$ is a dialgebra.
The identities that are necessary for embedding of $D$ into
$\mathrm{stl}(3,D)$ are declared to be the identities of an
alternative dialgebra.

The main purpose of this paper is to present a general scheme
how to deduce what is a variety of dialgebras similar to a
given variety of ordinary algebras
(associative, Lie, alternative, etc.).
This structure is motivated by the relation of dialgebras
and conformal algebras.

A conformal algebra \cite{K1} is a linear space $C$
(over a field of zero charachteristic) with one linear operation
$T: C\to C$ and a countable family of bilinear products
$(\cdot \oo{n}\cdot): C\times C \to C$, where $n$ ranges over the set
$\mathbb Z_+$ of non-negative
integers, such that
\begin{itemize}
\item
 for every $a,b\in C$ only a finite number of $a\oo{n}b$ is nonzero;
\item
 $Ta\oo{n} b = -n a\oo{n-1}b$, $a\oo{n}Tb = T(a\oo{n} b) + na \oo{n-1}b$.
\end{itemize}

Every conformal algebra can be constructed as a subspace of the
 formal power series space $A[[z,z^{-1}]]$
over an appropriate ordinary algebra $A$, with
respect to
\[
(Ta)(z) = \frac{d}{dz}a(z),\quad
(a\oo{n} b)(z)=\Ress_{w=0}a(w)b(z)(w-z)^n.
\]
Here $\Ress_{w=0} f(w,z)$ stands for the coefficient of $w^{-1}$
in $f(z,w)\in A[[z,z^{-1},w,w^{-1}]]$.

In \cite{Ro1} it was shown that for an arbitrary conformal algebra $C$
there exists unique ordinary algebra $\Coeff C$
(called coefficient algebra)
which is universal in the following sense:
$\Coeff C[[z,z^{-1}]]$ contains $C$, and
for any algebra $A$
such that $A[[z,z^{-1}]]$ contains $C$ as above there exists a
homomorphism $\Coeff C \to A$
such that the natural expansion $\Coeff C[[z,z^{-1}]]\to A[[z,z^{-1}]]$
acts on $C$ as an identity.

A conformal algebra $C$ is said to be associative (Lie, alternative,
etc.) \cite{Ro1}
if so is $\Coeff C$. A more general approach was developed in
\cite{Ko3}, see Section~\ref{sec4.1} for more details.
For example, $C$ is associative if
\[
a\oo{n}(b \oo{m} c) = \sum\limits_{s\ge 0} \binom{n}{s}
(a\oo{n-s} b) \oo{m+s} c,
\quad a,b,c\in C, \ n,m\in \mathbb Z_+.
\]

It is not difficult to note that an associative conformal algebra
$C$
endowed with two operations
\[
a\rbar b = a\oo{0} b,\quad
a\lbar b = \sum\limits_{s\ge 0}\frac{1}{s!}(-T)^{s}(a\oo{s}b)
\]
satisfies \eqref{assD1}--\eqref{assD3}, i.e., this
is an associative dialgebra. The same is true for alternative
conformal algebras.
However, it is not obvious why every associative (or alternative)
dialgebra can be obtained in this way. This is why we define what is a
variety of dialgebras using the language of operads, and then prove
that every dialgebra from a variety $\Var $ can be embedded
into a conformal algebra (more generally, pseudo-algebra)
from the variety $\Var $ of conformal algebras.

\section{Operads and algebras}

\subsection{Partitions}

Suppose $m\ge n\ge 1$ are two integers.
An ordered $n$-tuple of integers
$\pi=(m_1,\dots ,m_n)$, $m_i\ge 1$,
is called an $n$-{\it partition} of $m$
if  $m_1+\dots + m_n = m$.
The set of all such partitions is denoted by
$\Pi(m,n)$.

A partition $\pi\in \Pi(m,n)$
determines a one-to-one correspondence
between the sets
$\{1,\dots, m\}$ and $\{(i,j)\mid i=1,\dots, n, \,
 j=1,\dots, m_i\}$ by the rule
\[
 (i,j) \leftrightarrow k=(i,j)^\pi, \quad k=m_1+\dots + m_{i-1} + j.
\]
Given two partitions
$\tau=(p_1,\dots ,p_m)\in \Pi(p,m)$,
$\pi = (m_1,\dots ,m_n)\in \Pi(m,n)$,
define
$\tau\pi \in \Pi(p,n)$ in a natural way:
\[
  \tau\pi = (p_1+\dots+p_{m_1}, p_{m_1+1}+\dots + p_{m_1+m_2},
\dots, p_{m-m_n+1} +\dots+p_m).
\]
If $\tau\pi = (q_1,\dots, q_n)$
then
for any $i=1,\dots, n$ denote by
$\tau\pi_i$
the subpartition
$(p_{m_1+\dots+m_{i-1}+1}, \dots, p_{m_1+\dots+m_i}) \in \Pi(q_i,m_i)$.

We will also consider the (right) action of the symmetric
group $S_n$ on $\Pi(m,n)$:
if $\pi = (m_1,\dots, m_n)$, $\sigma \in S_n$
 then
 $\pi\sigma  = (m_{1\sigma^{-1}}, \dots, m_{n\sigma^{-1}})$.

\subsection{Operads}

Consider the following collection of data:
\begin{itemize}
\item[1)]
  a family of linear spaces
  $\{\mathcal C(n)\}_{n\ge 1}$;
\item[2)]
 a family of linear maps (called composition rule)
\[
  \Comp^\pi :\mathcal C(n)\otimes \mathcal C(m_1)\otimes \dots
  \otimes \mathcal C(m_n) \to \mathcal C(m),
\]
where $\pi=(m_1,\dots ,m_n) \in \Pi(m,n)$.
\end{itemize}

This collection is called an {\em operad\/}
(denoted by $\mathcal C$) if the following properties hold.
\begin{itemize}
\item[(O1)] The composition rule is associative, i.e.,
if
$\tau=(p_1,\dots ,p_m)\in \Pi(p,m)$,
$\pi = (m_1,\dots ,m_n)\in \Pi(m,n)$,
and
\[
\psi_j\in \mathcal C(p_j),
\
\chi_i \in \mathcal C(m_i),
\
\varphi \in \mathcal C(n), \quad
  i=1,\dots, n,
\ j=1,\dots, m,
\]
then
\begin{multline}
\label{P-assoc}
\Comp^\tau \big(\Comp^\pi \big(\varphi ,(\chi_i)_{i=1}^n\big),
 (\psi_j)_{j=1}^m \big) \\
=
\Comp^{\tau\pi}\big(\varphi ,
\big(\Comp^{\tau\pi_i} \big(\chi_i, (\psi_{(i,t)^\pi})_{t=1}^{m_i}
\big)_{i=1}^n\big)
\big).
\end{multline}

\item[(O2)]
There exists a ``unit'' $\idd \in \mathcal C(1)$
such that
\begin{equation}  \label{P-unit}
\Comp^{(1,\dots, 1)} (f, \idd ,\dots , \idd ) =
\Comp^{(n)} (\idd_A, f) = f
\end{equation}
for every $f\in \mathcal C(n)$.
\end{itemize}

An important example of an operad is provided by the following data.
Suppose $\mathcal C(n)= \Bbbk S_n$.
If $\sigma \in S_n$, $\pi=(m_1,\dots, m_n)\in \Pi(m,n)$,
$\tau_i \in S_{m_i}$,  $i=1,\dots, n$,
then define
\[
 \Sigma= \Comp^\pi(\sigma, \tau_1,\dots, \tau_n)\in S_m
\]
by the rule: if $k=(i,j)^\pi \in \{1,\dots, m\}$
then $k\Sigma = (i\sigma, j\tau_i)^{\pi\sigma }$.

The principle of constructing such compositions is clearly shown by the
following simple example:
\[
\Comp^{(3,2,4)}((123), (132), (12), (234))=
\begin{pmatrix} 1&2&3&4&5&6&7&8&9\\
7&5&6&9&8&1&3&4&2\end{pmatrix}.
\]
Graphically, this composition can be expressed by the diagram

\begin{center}
\includegraphics{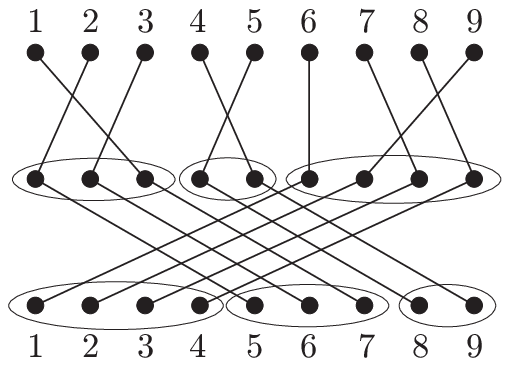}
\end{center}

The operad obtained is called the {\em operad of symmetries\/}
and denoted by $\Sym $.

Suppose $\mathcal C$ is an operad.
If $\mathcal C(n)$ is endowed with a right
action of the symmetric group $S_n$ for all $n\ge 1$,
 and for any
$\pi=(m_1,\dots, m_n) \in \Pi(m,n)$,
$\varphi \in \mathcal C(n)$,
$\sigma \in S_n$,
$\psi_i\in \mathcal C(m_i)$,
$\tau_i \in S_{m_i}$,
$i=1,\dots, n$,
we have
\[
 \Comp^{\pi\sigma } \big( \varphi^\sigma,
\psi_{1\sigma^{-1}}^{\tau_{1\sigma^{-1}}},
\dots,
\psi_{n\sigma^{-1}}^{\tau_{n\sigma^{-1}}}
\big)
=
\Comp^\pi (\varphi, \psi_1,\dots, \psi_n)
^{\Comp^\pi(\sigma, \tau_1,\dots, \tau_n)},
\]
then the operad $\mathcal C$ is said to be {\em symmetric}.

For example, $\Sym $ is a symmetric operad with respect to
the regular right action of $S_n$ on $\Sym (n)$:
$\varphi^\sigma = \varphi\sigma$, $\sigma\in S_n$,
$\varphi\in \Sym (n)=\Bbbk S_n$.

Let $\mathcal C_1$, $\mathcal C_2$ be two operads.
A {\em functor\/} $F$ from $\mathcal C_1$ to $\mathcal C_2$
is a family of linear maps
$F_n: \mathcal C_1(n) \to \mathcal C_2(n)$, $n\ge 1$,
which preserve the composition rule and the unit.
If $\mathcal C_i$ are symmetric operads, and
each $F_n$ is $S_n$-invariant then
$F$ is said to be symmetric.

If $\mathcal C_1$, $\mathcal C_2$ are two (symmetric)
operads, then the collection of spaces
$\mathcal C(n) = \mathcal C_1(n)\otimes \mathcal C_2(n)$,
$n\ge 1$,
defines a (symmetric) operad with respect to
the component\-wise composition (and symmetric group action).
We will denote this operad by $\mathcal C_1\otimes \mathcal C_2$.

\subsection{Examples of operads}

We have already considered the operad $\Sym $ built on symmetric groups.
Let us now state some well-known examples of operads that will be used later.
Throughout the rest of the subsection, $\pi$
stands for a partition $(m_1,\dots, m_n)\in \Pi(m,n)$,
$\sigma $ is a permutation from~$S_n$.

\subsubsection{Operad $\mathcal E$}
Set $\mathcal E(n) = \Bbbk^n$, and
let $\big\{e_1^{(n)}, \dots, e_n^{(n)}\big\}$
be the standard basis of $\Bbbk^n$.
Define a composition rule and $S_n$-action as follows:
\[
\Comp^{\pi}\big(
e_i^{(n)}, e_{j_1}^{(m_1)}, \dots, e_{j_n}^{(m_n)} \big)
= e^{(m)}_{m_1 + \dots + m_{i-1} + j_i},
\]
\[
\sigma : e_i^{(n)} \mapsto e_{i\sigma}^{(n)}.
\]
This is a symmetric operad.
In fact, $\mathcal E(n)$ can be considered as an image of
$\Sym(n)$ under the map $F_n:\Bbbk S_n\to \Bbbk^n$
given by $F_n(\sigma)= \sigma S_{n-1}$, assuming
$(in)S_{n-1}$ is identified with $e^{(n)}_i$.
The family of maps $F_n$, $n\ge 1$, defines
a functor from $\Sym $ to $\mathcal E$, however, non-symmetric.

\subsubsection{Operads $\Vec_{\Bbbk}$}
Let $V$ be a linear space over a field $\Bbbk$.
Define an operad which is also denoted by $V$ as follows:
\[
 V(n) = \Hom (V^{\otimes n}, V),\quad n\ge 1,
\]
$\Comp $ is the ordinary composition of multi-linear maps.
If $f\in V(n)$ then set
\[
 f^\sigma(v_1,\dots, v_n) = f(v_{1\sigma}, \dots, v_{n\sigma}),
\quad v_i\in V.
\]
This is a symmetric operad. We will denote the class of operads
obtained in this way by $\Vec_{\Bbbk }$.

\subsubsection{Operads $\Hmod$}
Suppose $H$ is a cocommutative bialgebra with
a coproduct $\Delta: H\to H\otimes H$.
We will use the standard Sweedler's notation:
\[
\Delta(h)=h_{(1)}\otimes h_{(2)},
\quad
\Delta^2(h)=(\Delta\otimes \idd_H)\Delta(h)=
(\idd_H\otimes \Delta)\Delta(h)= h_{(1)}\otimes h_{(2)}\otimes h_{(3)},
\]
etc.

For every left $H$-module $M$ one may consider
\begin{equation} \label{eq:Hmod}
  M(n) = \Hom_{H^{\otimes n}} (M^{\otimes n}, H^{\otimes n}\otimes _H M),
\end{equation}
where $H^{\otimes n}$ is considered as the outer product of
regular right $H$-modules.

An arbitrary
$f\in M(n)$ can be expanded to a map
\[
 f : (H^{\otimes m_1}\otimes_H M)\otimes \dots \otimes
 (H^{\otimes m_n}\otimes_H M) \to H^{\otimes m}\otimes_H M
\]
by the following rule:
if
$G_i\in H^{\otimes m_i}$, $a_i\in M$, $i=1,\dots, n$,
then
\begin{multline}\label{eq:*-exp}
f(G_1\otimes_H a_1,\dots, G_n\otimes_H a_n) \\
=
 (G_1\otimes \dots \otimes G_n\otimes _H 1)
 (\Delta^{m_1-1}\otimes \dots \otimes \Delta^{m_n-1}\otimes _H\idd_M)
 f(a_1,\dots, a_n).
\end{multline}
This expansion makes clear what is a composition of
maps \eqref{eq:Hmod}.
The action of $S_n$ on $M(n)$ can be defined by
\[
 \varphi^\sigma (a_1,\dots, a_n) = (\sigma\otimes _H \idd_M)
\varphi (a_{1\sigma},\dots, a_{n\sigma}),
\]
where $(h_1\otimes \dots \otimes h_n)^\sigma
 = (h_{1\sigma^{-1}}\otimes \dots \otimes h_{n\sigma^{-1}})\in H^{\otimes n}$.

We will denote the operad obtained in this way by the same symbol $M$.
The class of all such operads is denoted by $\Hmod$.

\subsubsection{Operad $\Alg$}
Let $X=\{x_1,x_2,\dots \}$ be a countable set of symbols.
Define $\Alg(n)$ to be the linear span
of all non-associative words obtained from
$x_1\dots x_n$ by some bracketing.
The natural composition rule
\begin{equation}\label{comp_Alg}
 \Comp^\pi (u, v_1,\dots, v_n)
=
u(v_1(x_{(1,1)^\pi}, \dots , x_{(1,m_1)^\pi}),
\dots ,
v_n(x_{(n,1)^\pi}, \dots , x_{(n,m_n)^\pi})
)
\end{equation}
defines a (non-symmetric) operad denoted by $\Alg$.

Consider $\Alg_S(n)=\Alg(n)\otimes \Bbbk S_n$ endowed with a composition
rule
\begin{multline}\label{comp_Alg_S}
 \Comp^\pi(u\otimes \sigma, v_1\otimes \tau_1,
\dots, v_n\otimes \tau_n) \\
=\Comp^{\pi\sigma^{-1} }(u, v_{1\sigma},\dots, v_{n\sigma})
 \otimes
\Comp^{\pi\sigma^{-1}} (\sigma, \tau_{1\sigma}, \dots, \tau_{n\sigma}).
\end{multline}
Hereafter, we will identify $(x_1\dots x_n)\otimes \sigma $
with the word $(x_{1\sigma}\dots x_{n\sigma })$.
Under the natural action of $S_n$, $\Alg_S$ is a symmetric operad.

If $\mathcal C$ is a symmetric operad and
$\Phi$ is a functor from $ \Alg$ to $\mathcal C$
then $\Phi $ can be extended to a functor
from $\Alg_S$ to $\mathcal C$ by the obvious
rule
\[
  \Phi_n(u\otimes \sigma ) = \Phi_n(u)^\sigma, \quad u\in \Alg(n), \
\sigma\in S_n.
\]

\subsubsection{Operad $\Var\Alg$}\label{sec:VarAlg}
Let $\Var $ be a homogeneous variety of (non-associative) algebras defined
by a family of polylinear identities.
Denote by $\Bbbk_{\Var}\{ X \} $ the free
algebra in this variety generated by $X=\{x_1,x_2,\dots \}$.
There exists an ideal $I_{\Var }$ of the free
algebra $\Bbbk \{X\}$ such that
$\Bbbk_{\Var}\{X\} \simeq \Bbbk\{X\}/I_{\Var}$.

Each space $\Alg_S(n) $, $n\ge 1$,
can be considered as a subspace of
$\Bbbk\{X\}$. Denote
\[
  \Var\Alg(n) = \Alg_S(n) / (\Alg_S(n)\cap I_{\Var}).
\]
Since the variety $\Var $ is homogeneous,
the composition map and the symmetric group
action are well-defined on the family $\{\Var\Alg(n)\}_{n\ge 1}$,
so we obtain a symmetric operad.
The family of projections $\Alg_S(n) \to \Var\Alg(n)$
is a functor from $\Alg_S$ to $\Var\Alg$, we will also
denote it by~$\Var $.

\subsection{Algebras and pseudo-algebras}\label{sec4.1}

Let $A$ stands for a (non-associative, in general)
algebra over~$\Bbbk $.
Then there exists a functor
(also denoted by $A$)
from $\Alg $ to the operad $A\in \Vec_{\Bbbk}$ built on the underlying
space.
Indeed, for every $u=(x_1\dots x_n)\in \Alg(n)$
define
$A_n(u): a_1\otimes \dots \otimes a_n \to (a_1\dots a_n)\in A$.

Conversely, every functor from $\Alg $ to $A\in \Vec_{\Bbbk}$
defines an operation $A\otimes A\to A$ which makes $A$ to be
an algebra.

Therefore, the notion of an algebra can be identified with
the notion of a functor from $\Alg $ to an operad from
$\Vec_{\Bbbk }$. Similarly, a pseudo-algebra
\cite{BDK}
is just a functor from $\Alg $ to an operad from
$\Hmod$. Since $\Vec_{\Bbbk}$ and $\Hmod $ consist of symmetric operads,
a (pseudo-)algebra can be considered as a functor defined
on $\Alg_S(n)$.

This is now clear how to define varieties of pseudo-algebras.
Indeed, suppose $\Var $ is a variety as in~\ref{sec:VarAlg}.
An (ordinary) algebra $A$ belongs to $\Var $
if and only if there exists a functor $\bar A:\Var\Alg \to A\in \Vec_{\Bbbk}$
such that
$A = \Var\circ \bar A $.
For pseudo-algebras, one should replace $\Vec_{\Bbbk}$ with
$\Hmod$ to obtain a definition of what does it mean
that a pseudo-algebra $C$ belongs to the variety $\Var $
\cite{Ko4}.

This approach leads to
pseudo-algebraic identities as in \cite{Ko3}.
If $\Var $ is defined by a family $\Sigma $ of
polylinear identities then a pseudo-algebra $C$
is said to be a $\Var$-pseudo-algebra if $C$
satisfies the identities
\[
  t^*(x_1,\dots, x_n), \quad t\in \Sigma,
\]
where $t^*$ is obtained from $t$ by the following way:
every monomial $\alpha(x_{1\sigma}\dots x_{n\sigma })$,
$\alpha\in \Bbbk$, $\sigma \in S_n$, that appears
in $t$ should be replaced with
$\alpha (\sigma \otimes_H \idd)(x_{1\sigma}*\dots *x_{n\sigma })$.
The expression $(x_{1\sigma }*\dots *x_{n\sigma})$ can
be evaluated at
$a_1,\dots, a_n\in C$ by \eqref{eq:*-exp},
its value lies in $H^{\otimes n}\otimes_H C$.


\section{Dialgebras and corresponding operads}

Consider an operad $\Dialg $, analogous to $\Alg $, generated by two binary
 products $\lbar $, $\rbar $.
The space $\Dialg(n)$ is the linear span of all terms obtained
from $x_1\dots x_n$ by some bracketing with operations~$\lbar $,~$\rbar $,
and the same rule as \eqref{comp_Alg} defines compositions on $\Dialg $.

By $\Dialg_S$ we denote the symmetric operad obtained from $\Dialg $
in the same way as $\Alg_S $ is obtained from $\Alg $. The operad
$\Dialg $ (resp., $\Dialg_S$) has the following property:
for any (symmetric) operad $\mathcal C$ and for any
$\varphi_1, \varphi_2 \in \mathcal C(2)$ there exists a
unique (symmetric) functor $\Psi : \Dialg \to \mathcal C$
($\Dialg_S \to \mathcal C$) such that
$\Psi_2(x_1\rbar x_2)=\varphi_1$,
$\Psi_2(x_1\lbar x_2)=\varphi_2$.
In particular, for $\mathcal C = \Alg_S\otimes \mathcal E$ there exists
$\Psi: \Dialg_S \to \Alg_S\otimes \mathcal E$
such that
\begin{equation}\label{eq:Psi}
\Psi_2(x_1\rbar x_2)=x_1x_2\otimes e_2^{(2)},
\quad
\Psi_2(x_1\lbar x_2)=x_1x_2\otimes e_1^{(2)}.
\end{equation}
It is clear that the functor $\Psi $ defined by \eqref{eq:Psi}
is full, i.e., $\Psi_n(\Dialg_S(n))=\Alg_S(n)\otimes \mathcal E(n)$.
Throughout the rest of the paper, $\Psi $ stands for the functor
$\Dialg_S \to \Alg_S\otimes \mathcal E$ defined by \eqref{eq:Psi}.

\begin{defn}
A dialgebra $A$ is said to be a 0-{\em dialgebra\/} if
$A$ satisfies the identities
\begin{equation}\label{eq:0-id}
(x\lbar y)\rbar z = (x\rbar y)\rbar z,
\quad
x\lbar (y \rbar z) = x\lbar (y \lbar z).
\end{equation}
\end{defn}

The purpose of this section is to prove the following statement.

\begin{thm}\label{thm2.1}
For any 0-dialgebra $A$,
considered as a functor $\Dialg_S \to A\in \Vec_{\Bbbk}$,
there exists a symmetric
functor $\bar A: \Alg_S\otimes \mathcal E \to A$
such that $A=\Psi\circ \bar A$.
\end{thm}

Let $H$ be a cocommutative Hopf algebra with a counit
$\varepsilon $ and an antipode~$S$.

If $C$ is a left $H$-module then $H^{\otimes n}\otimes_H C \simeq
H^{\otimes n-1}\otimes  C$ as a linear space \cite{BDK}.
We will denote the isomorphism by $i_n$ ($n\ge 1$), namely,
\[
i_n(h_1\otimes \dots \otimes h_{n} \otimes _H a)
=
h_1S(h_{n(1)})\otimes \dots \otimes h_{n-1}S(h_{n(n-1)})\otimes h_{n(n)}a,
\]
$h_i\in H$, $a\in C$.

For an operad $C\in \Hmod $ denote by
$C^{(0)}$ the operad from $\Vec_{\Bbbk}$ built on $C$
as on linear space.

\begin{lem}
A family of linear maps
\[
\varepsilon_n : C(n)\to C^{(0)}(n), \quad n\ge 1,
\]
given by
\[
\varepsilon_n(f)(a_1,\dots, a_n)
 = (\varepsilon^{\otimes n-1}\otimes \idd)i_n(f(a_1,\dots, a_n)),
 \quad
 f\in C(n),
\]
defines a (non-symmetric) functor from $C$ to $C^{(0)}$.
\end{lem}

\begin{proof}
The expansion rule \eqref{eq:*-exp}
immediately implies that
for all $f\in C(n)$, $g_i\in C(m_i)$, $i=1,\dots, n$,
$\pi=(m_1,\dots,m_n)\in \Pi(m,n)$, we have
\[
\varepsilon _m ( \Comp^{\pi}(f,g_1,\dots ,g_n))
=\Comp^\pi (\varepsilon _n(f), \varepsilon _{m_1}(g_1),
\dots ,\varepsilon _{m_n}(g_n)).
\]
Therefore, the family of maps $\varepsilon _n$, $n\ge 1$,
defines a functor (also denoted by~$\varepsilon $).
\end{proof}

Let $C$ stands for an $H$-pseudo-algebra considered as a functor
$\Alg_S \to C\in \Hmod$.
By the universal property of $\Dialg_S$ we have a uniquely
defined symmetric
functor from $\Dialg_S$ to $C^{(0)}\in \Vec_{\Bbbk}$
(also denoted by $C^{(0)}$)
such that
\begin{equation}\label{eq:C0}
\begin{gathered}
C^{(0)}_2 (x_1\rbar x_2) = \varepsilon_2(C_2(x_1x_2))\in C^{(0)}(2), \\
C^{(0)}_2 (x_1\lbar x_2) =
\varepsilon_2(C_2(x_2x_1))^{(12)} \in C^{(0)}(2).
\end{gathered}
\end{equation}
In particular, for every  $a,b\in C$
\begin{gather}
a\rbar b = C_2^{(0)} (x_1\rbar x_2)(a,b) =
(\varepsilon \otimes \idd)i_2 (a*b), \nonumber \\
a\lbar b=  C_2^{(0)} (x_1\lbar x_2)(a,b) = \varepsilon _2(C_2(x_2x_1))(b,a)
=(\varepsilon \otimes \idd)i_2\{a*b\},
\nonumber
\end{gather}
where $\{a*b\} = ((12)\otimes_H \idd)(a*b)$.

\begin{lem}
For every $H$-pseudo-algebra $C$ the corresponding
dialgebra $C^{(0)}$ given by \eqref{eq:C0} is a 0-dialgebra.
\end{lem}

\begin{proof}
It is straightforward to check that the identities \eqref{eq:0-id} hold
on $C^{(0)}$.
\end{proof}

\begin{prop}\label{prop:0-embed}
For every 0-dialgebra $A$ there exists an $H$-pseudo-algebra
$C=C(A)$ over $H=\Bbbk[T]$
such that $A $ is embedded into $C^{(0)}$.
\end{prop}

\begin{proof}
Here we consider $\Bbbk[T]$ with respect to the
usual Hopf algebra structure: $\Delta (T)=T\otimes 1 + 1\otimes T$,
$\varepsilon(T)=0$, $S(T)=-T$.

Let $A$ be a 0-dialgebra. Denote $\langle a,b\rangle = a\rbar b - a\lbar b$.
It follows from \eqref{eq:0-id} that
$\langle a,b\rangle \rbar c = 0$,
$a\lbar \langle b,c\rangle = 0$
for all $a,b,c \in A$.

Consider the linear space
\[
C = (H\otimes A)\oplus \left((A\otimes A)/\Span \{
\langle a,b \rangle \otimes \langle c,d\rangle \mid
a,b,c,d \in A
\}  \right).
\]
The first summand is a free $H$-module, the second one can be
endowed with an $H$-module structure as follows.
Define $T: A\otimes A \to A$ by
\begin{equation} \label{eq:Td}
T(a\otimes b) = \langle a,b\rangle,
\quad
a,b\in A,
\end{equation}
and note that $U=\Span \{
\langle a,b \rangle \otimes \langle c,d\rangle \mid
a,b,c,d \in A\} $ lies in the kernel of $T$.
Therefore, $T$ is well-defined on $(A\otimes A)/U$.

Later we will identify
$1\otimes a\in H\otimes A$ with $a$,
and write $a\otimes b$ ($a,b\in A$) for its image in $C$.

The pseudo-product on $C$ can be defined as follows:
\begin{gather}
a*b = 1\otimes 1\otimes _H (a\rbar b) - T\otimes 1\otimes_H (a\otimes b),
\label{eq:*1} \\
a*(b\otimes c)= 1\otimes 1\otimes_H a\otimes \langle b,c\rangle
= 1\otimes 1\otimes_H a\otimes T(b\otimes c),
\label{eq:*2} \\
(a\otimes b)* c  = - 1\otimes 1\otimes _H (\langle a,b \rangle \otimes c)
= - 1\otimes 1\otimes _H (T(a\otimes b) \otimes c),
\label{eq:*3}\\
(a\otimes b)*(c\otimes d) = 0.
\label{eq:*4}
\end{gather}
This operation is well-defined and $H^{\otimes 2}$-linear.
It is easy to see that $A$ is a subdialgebra of $C^{(0)}$.
\end{proof}

\begin{lem}\label{lem:l5}
Suppose $C$ is an $H$-pseudo-algebra.
Consider a family of linear maps
\[
C_n^{(\varepsilon)} : \Alg_S(n)\otimes \mathcal E(n) \to C^{(0)}(n),
\quad n\ge 1,
\]
defined by the following rule:
if
$C_n(f)(a_1,\dots, a_n) = \sum\limits_j h_{1j}\otimes \dots \otimes h_{nj}
 \otimes _H b_j$,
$f\in \Alg_S(n)$,
$a_i, b_j\in C$,
$h_{ij}\in H$,
then
\[
C_n^{(\varepsilon)}(f\otimes e_i)
 (a_1,\dots, a_n) =
 \sum\limits_j \varepsilon(h_{1j}\,\overset{i}{\widehat\dots}\, h_{nj})
 h_{ij}b_j.
\]
These maps define a symmetric functor
$C^{(\varepsilon)}: \Alg_S\otimes \mathcal E \to C^{(0)}\in \Vec_{\Bbbk}$.
\end{lem}

\begin{proof}
This is straightforward to check that
\[
  \big ( C_n^{(\varepsilon)} \big (f\otimes e_i^{(n)}\big ) \big)^\sigma
=
C_n^{(\varepsilon)} \big (f^\sigma \otimes e_{i\sigma }^{(n)} \big ),
\quad
f\in \Alg_S(n),\ i=1,\dots, n,\ \sigma \in S_n.
\]
It remains to
prove that $C_n^{(\varepsilon)}$ preserve compositions.
Since the composition rule is equivariant and associative, it is enough to
consider
$\pi =(n, m-n)\in \Pi(m,2)$,
$x_1x_2\in \Alg(2)$,
$u_1\in \Alg(n)$, $u_2\in \Alg(m-n)$, and
\begin{multline}
C_m^{(\varepsilon)}\big(
 \Comp^\pi \big (
 x_1x_2\otimes e_i^{(2)}, u_1\otimes e_{k}^{(n)},
  u_2\otimes e_{l}^{(m-n)} \big ) \big )   \\
=
C_m^{(\varepsilon)}\big ( f\otimes e_{p}^{(m)} \big )  ,
\quad
p = \begin{cases}
      k, & i=1, \\
      n+l, & i=2,
   \end{cases}
                                                \label{eq:e-comp}
\end{multline}
where $f=u_1(x_1,\dots, x_n)u_2(x_{n+1},\dots, x_m)$.
Now it is easy to note that
\[
 \Comp^{(n,m-n)} \big(C_2^{(\varepsilon)}
  \big (x_1x_2\otimes e_i^{(2)}\big ),
 C_n^{(\varepsilon)} \big (u_1\otimes e_{k}^{(n)} \big ),
 C_{m-n}^{(\varepsilon)} \big (u_2\otimes e_{l}^{(m-n)} \big )
  \big)
\]
coincides with \eqref{eq:e-comp}.
\end{proof}

\begin{proof}[Proof of Theorem~\ref{thm2.1}]
Let $A$ be a 0-dialgebra.
Denote by
$C$ its enveloping pseudo-algebra from Proposition~\ref{prop:0-embed}.
The functor $C: \Alg_S \to C\in \Hmod $
induces two functors
$C^{(0)}: \Dialg_S \to C^{(0)}$,
$C^{(\varepsilon)}: \Alg_S\otimes \mathcal E \to C^{(0)}$.
It is easy to see that
$\Psi \circ C^{(\varepsilon)} = C^{(0)}$.
Therefore, if $\Psi_n(f)=0$ for some $f\in \Dialg_S(n)$
then $C_{n}^{(0)}(f) = 0$ in $C^{(0)}(n)$. But $A$
is embedded into $C^{(0)}$, hence, $f$ is an identity
on $A_0$.

Therefore, a polynomial
$f\in \Dialg_S(n)$
is an identity on a 0-dialgebra $A$
provided that $\Psi_n(f)=0$.
Hence, there exists a family of linear maps
$\bar A_n: \Alg_S(n)\otimes \mathcal E(n) \to A(n)$
such that $A_n = \Psi_n\otimes \bar A_n$, $n\ge 1$.
This family defines a symmetric functor $\bar A$ from $\Alg_S\otimes \mathcal E$
to $A\in \Vec_{\Bbbk}$.
\end{proof}

One may construct an operad 0-$\Dialg$ related to the
variety of 0-dialgebras in the same way as $\Var\Alg $ relates
to a variety $\Var $.

\begin{cor}
The operad 0-$\Dialg $ is equivalent to $\Alg_S\otimes \mathcal E$.
\end{cor}

\section{Varieties of dialgebras}

Suppose $\Var $ is a homogeneous variety defined by a
family $\Sigma $ of polylinear identities
in $X=\{x_1,x_2,\dots  \}$. Then
$I_{\Var} $ is an ideal of $\Bbbk\{X\} $
generated by $\{t(f_1,\dots, f_n ) \mid t(x_1,\dots, x_n)\in \Sigma,\,
f_i\in \Bbbk\{X\}\}$.

Theorem~\ref{thm2.1} implies that every 0-dialgebra $A$
is actually a functor from $\Alg_S\otimes \mathcal E$
to $A\in \Vec_{\Bbbk}$. In view of this observation, the
following definition is quite natural.

\begin{defn}
A dialgebra $A:\Dialg_S \to A\in \Vec_{\Bbbk}$
 is said to be a member of the variety
$\Var $ of dialgebras
(or $\Var$-{\em dialgebra})
if there exists a symmetric functor
\[
\bar A: \Var\Alg \otimes \mathcal E \to A\in \Vec_{\Bbbk },
\]
such that $\Psi\circ (\Var\otimes \idd)\circ \bar A = A$,
i.e., the diagram of functors
\[
\begin{CD}
\Dialg_S @>A>> A\rlap{${}\in \Vec_\Bbbk$} \\
@V\Psi VV  @AA \bar{A} A \\
\Alg_S\otimes \mathcal E @>\Var\otimes \idd>>
 \Var\Alg\rlap{${}\otimes \mathcal E$}
\end{CD}
\]
is commutative.
\end{defn}

\begin{prop}\label{prop:p7}
A dialgebra $A$ belongs to the variety $\Var $
if and only if $A$ is a 0-dialgebra that
satisfies the identities $\Psi_n^{-1}\big (t\otimes e_i^{(n)}\big )$,
where $t\in \Sigma $, $n=\deg t$, $i=1,\dots, n$.
\end{prop}

\begin{proof}
It is enough to prove the ``if" part only.
Since $A$ is a 0-dialgebra,
it can be considered as a functor
$A : \Alg_S\otimes \mathcal E \to A\in \Vec_{\Bbbk}$.

Suppose
$f\in (\Alg_S\otimes \mathcal E)(n)$,
$f = \sum\limits_{i=1}^{n} f_i\otimes e_i^{(n)}$,
$f_i\in \Alg_S(n)$.
If $(\Var_n\otimes \idd)(f)=0$
then $\Var_n(f_i) = 0$ for all $i=1,\dots, n$.
Therefore, it needs to be shown that if
$\Var_n(g)=0$ for some $g\in \Alg_S(n)$
then $A_n\big (g\otimes e_i^{(n)}\big ) = 0$
for all $i=1,\dots, n$.
Since $\Ker \Var_n = \Alg_S(n)\cap I_{\Var}$,
the element $g$ can be presented as
\[
 g=\sum\limits_k \Comp^{\pi_k} (a_k, b_{k1}, \dots , b_{kl_k})^{\sigma_k},
\]
where
$a_k\in \Alg_S(l_k)$,
$b_{kj}\in \Alg_S(m_{kj})$,
$\pi_k = (m_{k1},\dots, m_{kl_k})\in \Pi(n, l_k)$,
$\sigma_k\in S_n$,
and
 for each $k $ there exists $j_k\in \{1,\dots, l_k\}$
such that
\[
  b_{kj_k} = \Comp^{\tau_{k}} (t_{k}, c_{k1},\dots, c_{km_k} ),
\quad t_k\in \Sigma, \quad c_{km}\in \Alg_S(q_m).
\]
For every $i\in \{1,\dots, n\}$ one may present $g\otimes e_i^{(n)}$
as a composition
in $\Alg_S\otimes \mathcal E$ with the same $a_k$, $b_{kj}$, $t_k$,
$c_{km}$ in the first tensor factor
(since the composition rule in $\Alg_S\otimes \mathcal E$
is componentwise).
Therefore,
$A_n\big(g\otimes e_i^{(n)}\big ) =0$
provided that
$A_{m_k}\big (t_k\otimes e^{(m_k)}_p \big)=0$
for every
$p=1,\dots, m_k$.
\end{proof}

To deduce defining identities in the signature of a dialgebra, the following
construction is useful.
Consider the functor $\alpha:\Dialg \to \Alg\otimes \Sym$
defined by $\alpha_2(x_1\rbar x_2) = x_1x_2\otimes \idd_2$,
$\alpha_2(x_1\lbar x_2)=x_1x_2\otimes (12)$.
Denote by $E$ the following functor from $\Sym $ to $\mathcal E$:
$E_n(\sigma ) = e^{(n)}_{n\sigma^{-1}}$.

\begin{lem}\label{lemA}
If $u\otimes \sigma \in  \Dialg_S(n)$,
$\alpha_n(u)=v\otimes \tau$,
then
\begin{equation}\label{eq:Psi-1}
\Psi_n(u\otimes \sigma ) = (v\otimes \sigma)\otimes E_n(\tau)^\sigma.
\end{equation}
\end{lem}

\begin{proof}
Denote the right-hand side of \eqref{eq:Psi-1} by $\Phi_n(u\otimes\sigma )$.
It is enough to check that

1) $\Phi_n( (u\otimes \sigma ))^{\sigma _1} =\Phi_n(u\otimes \sigma \sigma _1)$;

2) the maps $\Phi_n$, $n\ge 1$, preserve composition.

The first is obvious. The second needs to be checked
for $u=u_1\rbar u_2$ and $u=u_1\lbar u_2$, $\sigma =\idd_n$,
$u_i\in \Dialg(n_i)$, $n_1+n_2=n$.

If $u=u_1\rbar u_2 $,
$\alpha _{n_i}(u_i)=v_i\otimes \tau_i$,
 then $\alpha _n(u)=v_1v_2\otimes \Comp^{(n_1,n_2)}(\idd_2, \tau_1, \tau_2)$.
Hence,
\[
\Phi_n(u\otimes \idd_n)
= v_1v_2\otimes \idd_n \otimes e_{n_1+n_2\tau_2^{-1}}^{(n)}.
\]
On the other hand,
\begin{multline}
 \Comp^{(n_1,n_2)} \big (
x_1x_2\otimes e_2^{(2)},
\Phi_{n_1}(v_1\otimes \idd_{n_1}), \Phi_{n_2}(v_2\otimes \idd_{n_2})
\big ) \\
=
(v_1v_2\otimes \idd_n) \otimes \Comp^{(n_1,n_2)} \big(
 e_2^{(2)}, e^{(n_1)}_{n_1\tau_1^{-1}}, e^{(n_2)}_{n_2\tau_2^{-1}}
 \big)
=
 v_1v_2\otimes \idd_n \otimes e_{n_1+n_2\tau_2^{-1}}^{(n)}.
\nonumber
\end{multline}

If $u=u_1\lbar u_2 $,
$\alpha _{n_i}(u_i)=v_i\otimes \tau_i$,
 then $\alpha _n(u)=v_1v_2\otimes \Comp^{(n_1,n_2)}((12), \tau_1, \tau_2)$.
Hence,
\[
\Phi_n(u\otimes \idd_n)
= v_1v_2\otimes \idd_n \otimes e_{n_1\tau_1^{-1}}^{(n)}.
\]
On the other hand,
\begin{multline}
 \Comp^{(n_1,n_2)} \big (
 x_1x_2\otimes e_1^{(2)},
\Phi_{n_1}(v_1\otimes \idd_{n_1}), \Phi_{n_2}(v_2\otimes \idd_{n_2}) \big ) \\
=
(v_1v_2\otimes \idd_n) \otimes \Comp^{(n_1,n_2)} \big (e_1^{(2)},
 e^{(n_1)}_{n_1\tau_1^{-1}}, e^{(n_2)}_{n_2\tau_2^{-1}} \big )
=
 v_1v_2\otimes \idd_n \otimes e_{n_1\tau_1^{-1}}^{(n)}.
\nonumber
\end{multline}

Therefore, the sequence $\Phi_n$, $n\ge 1$, defines a symmetric
functor $\Phi $ from $\Dialg_S$ to $\Alg_S\otimes \mathcal E$ that coincides
with $\Psi $ on $\Dialg_S(2)$. Hence, $\Phi =\Psi$.
\end{proof}

In order to get
$\Psi_n^{-1}\big (v\otimes \sigma\otimes e_i^{(n)}\big)$
for some word $v\in \Alg(n)$, $\sigma \in S_n$, $i\in \{1,\dots, n\}$,
it is enough to put the signs $\lbar $, $\rbar $ of dialgebra operations
on the word $v$ with the same bracketing in such a way that
for the dialgebraic word $v_i^d\in \Dialg(n)$ obtained
one has
$\alpha_n\big (v_i^d\big )= v\otimes \tau$, where $i\sigma^{-1}=n\tau^{-1}$.
Then $v_i^d\otimes \sigma \in \Dialg_S(n)$ would be a preimage of
$v\otimes \sigma\otimes e_i^{(n)}$ with respect to~$\Psi_n$.

The method is quite clear. Suppose $v\otimes \sigma=(x_{1\sigma}\dots x_{n\sigma})
\in \Alg_S(n)$. Then
$v_i^d\otimes \sigma = (x_{1\sigma}\rbar \dots \rbar x_i
 \lbar \dots \lbar x_{n\sigma})$,
i.e., the dashes are always directed to the variable~$x_i$.

Let us consider some examples of varieties.

\subsection*{Associative dialgebras}
In this case $\Sigma =\{(x_1,x_2,x_3)\}$,
where $(x,y,z)=(xy)z-x(yz)$.
Lemma~\ref{lemA} implies that
\begin{equation}\label{diass}
\begin{gathered}
 \Psi_3^{-1}\big ((x_1,x_2,x_3) \otimes e_1^{(3)}\big )
 = (x_1,x_2,x_3)_{\lbar}, \\
 \Psi_3^{-1}\big ((x_1,x_2,x_3) \otimes e_2^{(3)}\big )
 = (x_1,x_2,x_3)_{\times}, \\
 \Psi_3^{-1}\big ((x_1,x_2,x_3) \otimes e_3^{(3)}\big )
 = (x_1,x_2,x_3)_{\rbar},
\end{gathered}
\end{equation}
where
\begin{gather}\nonumber
(x,y,z)_{\lbar}=(x\lbar y)\lbar z-x\lbar (y\lbar z),
\quad
(x,y,z)_{\times}=(x\rbar y)\lbar z-x\rbar (y\lbar z), \\
(x,y,z)_{\rbar}=(x\rbar y)\rbar z-x\rbar (y\rbar z).
\nonumber
\end{gather}
Together with \eqref{eq:0-id}, the identities \eqref{diass}
are exactly those introduced in \cite{L1, L2}.

\subsection*{Commutative dialgebras}
Let $\Sigma =\{(x_1,x_2,x_3), \, x_1x_2-x_2x_1 \}$.
The first identity implies \eqref{diass} as before, the last
one is equivalent to
\begin{equation}\label{dicomm}
x_1\rbar x_2 - x_2\lbar x_1.
\end{equation}
Therefore, a commutative dialgebra $A$ can be considered as an ordinary algebra
with respect to $ab=a\rbar b$, $a,b\in A$.
In this case, we obtain an associative algebra that satisfies
\begin{equation}\label{Pcomm}
 [x_1,x_2]x_3 = 0.
\end{equation}
Conversely, an arbitrary associative  algebra
that satisfies \eqref{Pcomm} is actually a commutative dialgebra.

\subsection*{Alternative dialgebras}
Suppose
$\Sigma =\{(x_1,x_2,x_3)+ (x_2,x_1,x_3), \, (x_1,x_2,x_3)+(x_1,x_3,x_2)\}$.
We have already seen what is obtained from $(x_1,x_2,x_3)\otimes e_i^{(3)}$
under $\Psi_3^{-1}$.
Therefore, the identities $\Sigma $ give rise to
\begin{equation}\label{dialt}
\begin{gathered}
(x_1,x_2,x_3)_{\lbar } + (x_1,x_2,x_3)_{\times }^{(12)},\quad
(x_1,x_2,x_3)_{\rbar } + (x_1,x_2,x_3)_{\rbar }^{(12)},\\
(x_1,x_2,x_3)_{\lbar } + (x_1,x_2,x_3)_{\lbar }^{(23)},\quad
(x_1,x_2,x_3)_{\times } + (x_1,x_2,x_3)_{\rbar }^{(23)}.
\end{gathered}
\end{equation}
These relations are equivalent to those introduced in \cite{DL}.

\subsection*{Lie dialgebras}
If $\Sigma  = \{(x_1,x_2,x_3)- x_2(x_1x_3),\, x_1x_2+x_2x_1\}$
then the corresponding dialgebra identities
include
\[
  x_1\lbar x_2 + x_2\rbar x_1.
\]
A Lie dialgebra $A$ considered as an ordinary algebra with respect to
$[ab]= a\rbar b$,    $a,b\in A$, is just a left Leibniz algebra.
Conversely, every left Leibniz algebra $L$ is a Lie dialgebra with
respect to $a\rbar b = [ab]$, $a\lbar b = -[ba]$.
Therefore, a Lie dialgebra is just the same as a Leibniz algebra.

\subsection*{Jordan dialgebras}

Let $\Sigma $ consists of the commutativity identity and
\begin{multline}
x_1 (x_2 (x_3 x_4)) + (x_2 (x_1 x_3)) x_4
 + x_3 (x_2 (x_1 x_4))  \\
 = (x_1 x_2)(x_3 x_4) + (x_1 x_3)(x_2 x_4)
  + (x_3 x_2)  (x_1 x_4).
                               \label{J2.6}
\end{multline}
The last one is a polylinear Jordan identity.

Therefore, a Jordan dialgebra satisfies commutativity
relation $x_1\lbar x_2 = x_2\rbar x_1$
and the following four identities:
\begin{multline}
  x_1\lbar (x_2\lbar (x_3\lbar x_4))
+ (x_2\rbar (x_1\lbar x_3))\lbar x_4
+ x_3 \rbar(x_2\rbar (x_1\lbar x_4))   \\
 = (x_1\lbar x_2)\lbar(x_3\lbar x_4)
 + (x_1\lbar x_3)\lbar(x_2\lbar x_4)
 + (x_3 \rbar x_2)\rbar(x_1\lbar x_4) .\nonumber
\end{multline}
\begin{multline}
  x_1 \rbar(x_2 \lbar(x_3 \lbar x_4))
+ (x_2\lbar (x_1\lbar x_3))\lbar x_4
+ x_3 \rbar(x_2 \lbar(x_1\lbar x_4))   \\
 = (x_1\rbar x_2)\lbar(x_3 \lbar x_4)
 + (x_1\rbar x_3)\rbar(x_2\lbar x_4)
 + (x_3\rbar x_2)\lbar (x_1\lbar x_4) .\nonumber
\end{multline}
\begin{multline}
  x_1\rbar (x_2\rbar (x_3\lbar x_4))
+ (x_2\rbar (x_1\rbar x_3))\lbar x_4
+ x_3\lbar (x_2\lbar (x_1\lbar x_4))   \\
 = (x_1\rbar x_2)\rbar(x_3\lbar x_4)
 + (x_1\rbar x_3)\lbar(x_2\lbar x_4)
 + (x_3\lbar x_2)\lbar(x_1\lbar x_4) .\nonumber
\end{multline}
\begin{multline}
  x_1 \rbar(x_2\rbar (x_3\rbar x_4))
+ (x_2\rbar (x_1\rbar x_3))\rbar x_4
+ x_3\rbar (x_2\rbar (x_1\rbar x_4))   \\
 = (x_1\rbar x_2)\rbar(x_3\rbar x_4)
 + (x_1\rbar x_3)\rbar(x_2\rbar x_4)
 + (x_3\rbar x_2)\rbar(x_1\rbar x_4) .\nonumber
\end{multline}

Being expressed in terms of single operation, say,
$x y = x\rbar y$,
these identities together with
 \eqref{eq:0-id}
turn into
\begin{equation}
(x_1x_2)x_3 = (x_2x_1)x_3. \label{JL0}
\end{equation}
\begin{multline}
  ((x_4 x_3)x_2)x_1
+ x_4(x_2(x_3 x_1))
+ x_3(x_2(x_4 x_1))   \\
 = (x_4 x_3)(x_2 x_1)
 + (x_4 x_2)(x_3 x_1)
 + (x_3 x_2)(x_4 x_1) .\label{JL1}
\end{multline}
\begin{multline}
  x_1((x_4x_3)x_2)
+ x_4((x_3 x_1)x_2)
+ x_3((x_4x_1)x_2)   \\
 = (x_4 x_3)(x_1x_2)
 + (x_1x_3)(x_4 x_2)
 + (x_4 x_1)(x_3 x_2) .\label{JL2}
\end{multline}
Identities \eqref{JL1}, \eqref{JL2} and \eqref{JL0} define a variety
of algebras
that relate to Jordan algebras as Leibniz algebras
relate to Lie algebras. These algebras inherit the important property
of Jordan algebras: the commutator of two operators of
left multiplication
is a derivation of the algebra.

\section{Relation to pseudo-algebras}

Suppose $\Var $ is a homogeneous
 variety defined by a family
of polylinear identities~$\Sigma $.

\begin{thm}\label{thmVar0}
If $C$ is an $\Var$-pseudo-algebra
over a cocommutative Hopf algebra $H$
then $C^{(0)}$ is a $\Var $-dialgebra.
\end{thm}

\begin{proof}
Every identity  $t\in \Sigma $
of degree $n$ can be considered as an element
of $\Alg_S(n)$, i.e.,
$t=\sum\limits_{\sigma \in S_n}t_\sigma \otimes \sigma $,
$t_\sigma \in \Alg(n)$.

Let us fix a linear basis $\{h_i\}$ of $H$ and
some elements $a_1,\dots ,a_n\in C$. Denote
\[
C_n(t_\sigma\otimes \sigma)(a_1,\dots ,a_n) = \sum\limits_{i_1,\dots, i_{n-1}}
h_{i_1}\otimes \dots \otimes h_{i_{n-1}}\otimes 1 \otimes_H
 b_{i_1,\dots ,i_{n-1}; \sigma }.
\]
Since $C$ is a $\Var$-pseudo-algebra, we have
$\sum\limits_{\sigma \in S_n} b_{i_1,\dots ,i_{n-1};\sigma }=0$
for all $i_1,\dots ,i_{n-1}$.

By Lemma~\ref{lem:l5}
\[
C_n^{(\varepsilon )}\big (t_\sigma \otimes \sigma \otimes e_j^{(n)}\big )
(a_1,\dots, a_n) =
\sum\limits_{i_1,\dots ,i_{n-1}}
\varepsilon (h_{i_1} \,\overset{j}{\widehat\dots}\, h_{i_{n-1}})h_{i_j}
 b_{i_1,\dots ,i_{n-1};\sigma }
\]
for $j\ne n$, and
\[
C_n^{(\varepsilon )}(t_\sigma \otimes \sigma \otimes e_n^{(n)})
(a_1,\dots, a_n) =
\sum\limits_{i_1,\dots ,i_{n-1}}
\varepsilon (h_{i_1} \dots h_{i_{n-1}})
 b_{i_1,\dots ,i_{n-1};\sigma }.
\]
Therefore,
$C^{(0)}_n \big (\Psi^{-1}_n\big (t\otimes e_j^{(n)}\big )\big )
= C_n^{(\varepsilon )}\big (t\otimes e_j^{(n)} \big )=0$ in $C^{(0)}(n)$,
so, $C^{(0)}$ is a $\Var$-dialgebra.
\end{proof}

It turns out that the converse is also true: a $\Var$-dialgebra
can be embedded into a pseudo-algebra (over $H=\Bbbk[T]$)
that
belongs to the variety $\Var $ of pseudo-algebras.

Let $A$ be a $\Var $-dialgebra.
Since $A$ is in particular a 0-dialgebra, there exists
a pseudo-algebra $C$ such that $A$
is a subdialgebra of $C^{(0)}$.
The purpose of this section is to show that
one may choose $C$ from the variety $\Var $ of
pseudo-algebras.

Recall the construction from Proposition~\ref{prop:0-embed}:
consider
$C=C_{(0)}\oplus C_{(1)}$, where $C_0=H\otimes A$, $C_1=(A\otimes A)/U$,
$H=\Bbbk[T]$, $U=\Span \{\langle a,b\rangle\otimes \langle c,d\rangle
\mid a,b,c,d\in A\}$, and the pseudo-product is given by
\eqref{eq:*1}--\eqref{eq:*4}.

As before, we identify
$A$ and $C$ with the symmetric functors
$A:\Alg_S\otimes \mathcal E \to A\in \Vec_{\Bbbk}$,
$C:\Alg_S\to C\in \Hmod$, respectively.

We will denote by $T_i^{(n)}$ the element
$(T\otimes 1 \otimes \dots \otimes 1)^{(1i)}\in H^{\otimes n}$,
$i=1,\dots, n$, and let $T_0^{(n)}$ stands for $1\otimes \dots \otimes 1
\in H^{\otimes n}$.

\begin{lem}\label{lem:t-alg}
Suppose $t\in \Alg(n)$, $a_1,\dots, a_n\in A$.
Then
\[
C_n(t)(a_1, \dots, a_n)  =T_0^{(n)}\otimes_H x_0
-
\sum\limits_{i=1}^{n-1} T_i^{(n)}\otimes _H x_i,
\]
where
$x_0\in C_{(0)}$, $x_i\in C_{(1)}$,
\begin{gather}
x_0 = A_n\big (t\otimes e_n^{(n)}\big )(a_1,\dots, a_n), \label{eq:t-alg1}\\
Tx_i = A_n\big (t\otimes e_n^{(n)} - t\otimes e_i^{(n)}\big )(a_1,\dots, a_n),
\quad i=1,\dots, n-1.
   				\label{eq:t-alg2}
\end{gather}
\end{lem}

\begin{proof}
For $n=2$ the statement is clear from \eqref{eq:*1}.
Assume we have the required relations for all terms
of degree $k$, $k<n$.
Then  an arbitrary $t\in \Alg(n)$
can be presented as
$t=\Comp^{(m,n-m)}(x_1x_2, t_1, t_2)$,
$t_1\in \Alg(m)$, $t_2\in \Alg(n-m)$.
By the assumption,
\[
C_m(t_1)(a_1,\dots, a_m)
 = T_0^{(m)}\otimes_H x_0 - \sum\limits_{i=1}^{m-1}T_i^{(m)}\otimes_H x_i,
\]
\[
C_{n-m}(t_2)(a_{m+1},\dots, a_n)
 = T_0^{(n-m)}\otimes_H y_0 -
 \sum\limits_{j=1}^{n-m-1}T_j^{(n-m)}\otimes_H y_j.
\]
Direct computation shows that
\begin{multline}
C_n(t)(a_1,\dots, a_n)
=C_m(t_1)(a_1,\dots, a_m)*C_{n-m}(t_2)(a_{m+1},\dots, a_n) \\
=
T_0^{(n)}\otimes_H (x_0\rbar y_0) - \sum\limits_{i=1}^{m} T_i^{(n)}
\otimes_H (x_0\otimes y_0)  \\
- \sum\limits_{i=1}^{m-1} T_i^{(n)}\otimes_H
(-Tx_i)\otimes y_0 - \sum\limits_{j=1}^{n-m-1} T_{m+j}^{(n)}\otimes_H
 x_0\otimes Ty_j.
 \nonumber
\end{multline}
Therefore,
\[
C_n(t)(a_1,\dots, a_n) =T_0^{(n)}\otimes_H z_0 - \sum\limits_{i=1}^{n-1}
 T_i^{(n)}\otimes_H z_i,
\]
where
\[
z_i =
\begin{cases}
 x_0\rbar y_0, & i=0, \\
 x_0\otimes y_0 - Tx_i\otimes y_0, & i=1,\dots, m-1, \\
 x_0\otimes y_0, & i=m, \\
 x_0\otimes Ty_{i-m}, & i=m+1,\dots, n-1.
\end{cases}
\]
By the assumption,
\[
x_0 = A_m\big (t_1\otimes e_m^{(m)}\big )(a_1,\dots, a_m),
\quad
y_0 = A_{n-m}\big (t_2\otimes e_{n-m}^{(n-m)}\big )(a_{m+1}, \dots, a_n).
\]
Hence,
\begin{multline}
z_0 = A_n\big(\Comp^{\pi}\big (x_1x_2\otimes e_2^{(2)}, t_1\otimes e_m^{(m)},
t_2\otimes e_{n-m}^{(n-m)} \big )\big ) (a_1,\dots, a_n) \\
=A_n\big (t\otimes e_n^{(n)}\big )(a_1,\dots, a_n),
\nonumber
\end{multline}
where $\pi=(m,n-m)\in \Pi(n,2)$.
Similarly, for $i=1,\dots, m-1$
one has
\begin{multline}
Tz_i = x_0\rbar y_0 - x_0\lbar y_0 -Tx_i\rbar y_0 + Tx_i\lbar y_0 \\
=
A_n\big(\Comp^{\pi}\big (x_1x_2\otimes \big (e_2^{(2)}-e_1^{(2)}\big ),
t_1\otimes e_m^{(m)}, t_2\otimes e_{n-m}^{(n-m)}
\big )\big)(a_1,\dots, a_n) \\
-
A_n\big(\Comp^{\pi }\big (x_1x_2\otimes \big (e_2^{(2)}-e_1^{(2)}\big ),
t_1\otimes \big (e_m^{(n)}- e_i^{(m)}\big ),
t_2\otimes e_{n-m}^{(n-m)} \big )\big)(a_1,\dots, a_n) \\
=
A_n\big(t\otimes e_n^{(n)} - t\otimes e_m^{(n)} -
t\otimes e_n^{(n)} + t\otimes e_n^{(n)} + t\otimes e_m^{(n)}
-t\otimes e_i^{(n)} \big )(a_1,\dots, a_n) \\
=
A_n\big(t\otimes \big (e_n^{(n)}- e_i^{(n)}\big ) \big)(a_1,\dots, a_n).
\nonumber
\end{multline}
In the same way, one may proceed in the remaining cases
for $i=m, m+1,\dots, n$.
\end{proof}

\begin{lem}\label{lem:S-alg}
Suppose $t\in \Alg_S(n)$,
$t=u\otimes \sigma$, $u\in \Alg(n)$, $\sigma \in S_n$,
and let $a_1,\dots, a_n \in A$. Then
\[
C_n(t)(a_1,\dots, a_n)=
T_0^{(n)}\otimes _H x_{0,\sigma}
-\sum\limits_{i=1}^{n-1} T_i^{(n)}\otimes _H x_{i,\sigma},
\]
where
$x_0\in C_{(0)}$, $x_i\in C_{(1)}$,
\begin{gather}
x_{0,\sigma} = A_n\big (t\otimes e_{n\sigma^{-1}}^{(n)}\big )(a_1,\dots, a_n),
\label{eq:S-alg1} \\
x_{i,\sigma } = A_n\big (t\otimes \big (e_{n\sigma^{-1}}^{(n)}
 - e_{i\sigma^{-1}}^{(n)}\big )\big )(a_1,\dots, a_n),
 \quad i=1,\dots, n-1.
\label{eq:S-alg2}
\end{gather}
\end{lem}

\begin{proof}
By the definition of the operad $C\in \Hmod$,
\[
C_n(t)(a_1,\dots, a_n) = (\sigma\otimes_H\idd)
  C_n(u)(a_{1\sigma}, \dots, a_{n\sigma}).
\]
By Lemma~\ref{lem:t-alg},
$C_n(u)(a_{1\sigma}, \dots, a_{n\sigma })=
T_0^{(n)}\otimes_H y_0 - \sum\limits_{i=1}^{n-1} T_i^{(n)}\otimes _H y_i$,
where $y_0, y_1, \dots, y_n$
satisfy \eqref{eq:t-alg1}, \eqref{eq:t-alg2}.

If $n\sigma =n$ then the statement is obvious.
Suppose $n\sigma \ne n$.
Since $\sigma:T_i^{(n)} \mapsto T_{i\sigma}^{(n)}$, we have
\begin{multline}
C_n(t)(a_1,\dots, a_n)
=
T_0^{(n)}\otimes_H y_0
- \sum\limits_{i=1}^{n-1} T_{i\sigma}^{(n)}\otimes _H y_i\\
=
T_0^{(n)}\otimes_H y_0
- T_n^{(n)}\otimes_H y_{n\sigma^{-1}}
- \sum\limits_{\hbox{\scriptsize
 $\array{c} 1\le j\le n-1 \\ j\ne n\sigma\endarray$}}
  T_{j}^{(n)}\otimes _H y_{j\sigma^{-1}} \\
=
T_0^{(n)}\otimes_H (y_0-Ty_{n\sigma^{-1}})
+
\sum\limits_{i=1}^{n-1} T_i^{(n)}\otimes _H
 (y_{n\sigma^{-1}} -(1-\delta_{i,n\sigma}) y_{i\sigma^{-1}}).
 \nonumber
\end{multline}
Therefore,
\[
x_i = \begin{cases}
 y_0-Ty_{n\sigma^{-1}}, & i=0, \\
 - y_{n\sigma^{-1}} + y_{i\sigma^{-1}} , &  i\ne n\sigma, \\
 - y_{n\sigma^{-1}}, & i=n\sigma.
\end{cases}
\]
It remains to apply
\eqref{eq:t-alg1} and \eqref{eq:t-alg2} to obtain
\eqref{eq:S-alg1} and \eqref{eq:S-alg2}.
\end{proof}

\begin{lem}\label{lem:D-alg}
Suppose $t\in \Alg(n)$, $a_1, \dots, a_n\in A$, $b\in A$.
Then for every $i=1,\dots, n$
\[
C_n(t)(a_1,\dots, a_{i-1}, a_i\otimes b, a_{i+1}, \dots, a_n)
= T_0^{(n)} \otimes_H x_i,
\]
where $x_i\in C_{(1)}$,
\begin{gather}
Tx_i = A_{n+1}\big(f_i\otimes \big (e_{i+1}^{(n+1)} - e_i^{(n+1)}\big ) \big)
 (a_1,\dots, a_i, b, a_{i+1}, \dots, a_n),
 \label{eq:D-alg1} \\
f_i = t(x_1,\dots, x_{i-1}, x_ix_{i+1}, x_{i+2}, \dots, x_{n+1}).
\nonumber
\end{gather}
\end{lem}

\begin{proof}
Suppose $t\in \Alg(n)$. If $n=1$ then the statement is clear.
Assume we are done for any $k<n$. It is sufficient to show the required
relations for $t=(x_1\dots x_n)\in \Alg(n)$.
Let us present $t$ as
$t=\Comp^{(m,n-m)}(x_1x_2, t_1,t_2)$,
$1\le m\le n$,
$t_1=(x_1\dots x_m)_1$,
$t_2=(x_1\dots x_{n-m})_2$,
where $(\dots )_1$, $(\dots )_2$ are some bracketings.
Suppose $i\in \{1,\dots ,m\}$.
Then
\begin{multline}\nonumber
C_n(t)(a_1,\dots ,a_i\otimes b_i, \dots a_n) \\
=
 C_{m}(t_1)(a_1,\dots ,a_i\otimes b_i, \dots a_m)*
 C_ {n-m}(t_2)(a_{m+1},\dots , a_n).
\end{multline}
By the assumption
\[
C_m(t_1)(a_1,\dots, a_i\otimes b, \dots, a_m)
= T_0^{(m)} \otimes_H x_i,
\]
where $x_i\in C_{(1)}$  satisfies \eqref{eq:D-alg1},
and by Lemma~\ref{lem:t-alg}
\[
C_{n-m}(t_2)(a_{m+1}, \dots, a_n)  =T_0^{(n-m)}\otimes_H y_0
-
\sum\limits_{j=1}^{n-m-1} T_j^{(n-m)}\otimes _H y_j,
\]
$y_0=A_{n-m}\big (t_2\otimes e_{n-m}^{(n-m)} \big  )(a_{m+1},\dots ,a_n)$,
$y_j\in C_{(1)}$ for $j=1,\dots ,n-m$.

Since $C_{(1)}*C_{(1)} =0$, we have
\[
C_n(t)(a_1,\dots ,a_i\otimes b, \dots , a_n)
= \big (T_0^{(m)}\otimes_H x_i\big ) * \big (T_0^{(n-m)}\otimes _H y_0\big )
=-T_0^{(n)}\otimes_H (Tx_i\otimes y_0).
\]
Hence, it is enough to consider
$-T(Tx_i\otimes y_0)= Tx_i\lbar y_0 - Tx_i\rbar y_0$.
Recall that
\[
  Tx_i = A_{m+1}\big (f_i\otimes \big (e^{(m+1)}_{i+1} - e^{(m+1)}_i\big ) \big )
   (a_1,\dots ,a_i,b,\dots , a_m),
\]
where $f_i=t_1(x_1,\dots ,x_ix_{i+1},\dots ,x_{m+1})$.
Therefore,
\[
Tx_i\lbar y_0 - Tx_i\rbar y_0 = F_i(a_1,\dots ,a_i,b, \dots ,a_n),
\]
where
\begin{multline}
F_i=
\Comp^{(m+1,n-m)}\big(x_1x_2\otimes \big (e_1^{(2)} - e_2^{(2)}\big ),
 f_i\otimes \big (e^{(m+1)}_{i+1} - e^{(m+1)}_i\big ),
 t_2\otimes e_{n-m}^{(n-m)}\big)
\\
=
t_1(x_1,\dots ,x_ix_{i+1}, \dots x_{m+1})t_{2}(x_{m+2},\dots, x_n)
\otimes \big (e_i^{(n+1)}-e_i^{(n+1)}\big ) \\
=
t(x_1,\dots ,x_ix_{i+1}, \dots  x_{n+1})\otimes \big (e_i^{(n+1)}-e_i^{(n+1)}\big ).
\nonumber
\end{multline}

The remaining case $i\in \{m+1,\dots ,n\}$
is completely analogous.
\end{proof}

\begin{lem}\label{lem:DS-alg}
Suppose $t\in \Alg_S(n)$,
$t=u\otimes \sigma$,
$u\in \Alg(n)$, $\sigma \in S_n$, and let
 $a_1, \dots, a_n\in A$, $b\in A$.
Then for every $i=1,\dots, n$
\[
C_n(t)(a_1,\dots, a_{i-1}, a_i\otimes b, a_{i+1}, \dots, a_n)
= T_0^{(n)} \otimes_H x_i,
\]
where $x_i\in C_{(1)}$,
\begin{gather}
Tx_i =
 A_{n+1}\big(f_i\otimes \big (e_{i\sigma^{-1}+1}^{(n+1)}
  - e_{i\sigma^{-1}}^{(n+1)}\big ) \big)
 (a_1,\dots, a_i, b, a_{i+1}, \dots, a_n),
 \label{eq:D-alg2} \\
f_i = t(x_1,\dots, x_{i-1}, x_ix_{i+1}, x_{i+2}, \dots, x_{n+1}).
\nonumber
\end{gather}
\end{lem}

\begin{proof}
By the definition of the operad $C\in \Hmod$,
\[
C_n(t)(a_1,\dots, a_{i-1}, a_i\otimes b, a_{i+1}, \dots, a_n)
=
(\sigma \otimes_H \idd)C_n(u)(z_1, \dots, z_n),
\]
where $z_j = a_{j\sigma}$ for $j\ne i\sigma^{-1}$,
and
$z_{i\sigma^{-1}}=a_i\otimes b$.
By Lemma~\ref{lem:D-alg} we get the required relation.
\end{proof}

\begin{thm}\label{thm:t2}
For every  $\Var $-dialgebra $A$
there exists a $\Var $-pseudo-algebra $C_{\Var} (A)$
such that $A\subseteq C_{\Var}(A)^{(0)}$.
\end{thm}

\begin{proof}
Consider the pseudo-algebra $C=C(A)$ from
Proposition~\ref{prop:0-embed}. It is enough to show that
there exists an ideal $I$ of $C$ such that
$I\cap C_{(0)} =0$ and $C/I$ is a $\Var $-pseudo-algebra.

Recall that $\Sigma $ is the system of polylinear  defining
identities of the variety $\Var$.
For any $t\in \Sigma $ of degree $n\ge 2$,
consider the expressions
\begin{gather}
C_n(t)(a_1,\dots, a_n) = T_0^{(n)}\otimes_H x_0
 - \sum\limits_{i=1}^{n-1}T_i^{(n)}\otimes _H x_i,
                 \label{eq:C00}
\\
C_n(t)(a_1,\dots, a_{j-1}, a_j\otimes b, a_{j+1}, \dots, a_n)
=T_0^{(n)}\otimes_H y_j,
                 \label{eq:C1}
\end{gather}
for all $a_1,\dots, a_n, b\in A$, $j=1,\dots, n$.
By \eqref{eq:t-alg1} and \eqref{eq:S-alg1} we have $x_0=0$.

Denote by $I$ the linear span of all elements
 $x_i, y_j\in C_{(1)}$ from \eqref{eq:C00}, \eqref{eq:C1},
$i,j\ge 1$, for all $a_1,\dots, a_n,b\in A$
and for all $t\in \Sigma $.
It follows from Lemmas \ref{lem:t-alg}--\ref{lem:DS-alg}
and \eqref{eq:*2}, \eqref{eq:*3} that $I$ is actually an ideal of
the pseudo-algebra $C$.
Moreover, it is clear that $C_{\Var}(A)=C/I$
(as a pseudo-algebra) satisfies all of the identities
$t\in \Sigma $:  if two or more arguments
among $x_1,\dots, x_n\in C$ belong to $C_{(1)}$
then $C_n(t)(x_1,\dots, x_n)=0$ by \eqref{eq:*4}.
\end{proof}

\begin{rem}
Theorem~\ref{thm:t2} provides another proof of Proposition~\ref{prop:p7}.
\end{rem}

Indeed, if $A$ is a 0-dialgebra that satisfies
$t\otimes e_i^{(n)}$, $t\in \Sigma$, $\deg t=n$, $i=1,\dots, n$,
then there exists a $\Var $-pseudo-algebra $C$ such that $A\subset C^{(0)}$.
If $g\in \Alg_S(n)$, $\Var_n(g)=0$, then by Lemma~\ref{lem:l5}
we have $C_n^{(\varepsilon)}\big (g\otimes e_i^{(n)}\big )=0$
for all $i=1,\dots, n$.
Therefore, $g\otimes e_i^{(n)}$ is an identity on $C^{(0)}\supset A$.

\begin{prop}\label{prop:p16}
Let $A$ be a $\Var $-dialgebra, and let $P$ be a
$\Var$-pseudo-algebra over $H=\Bbbk[T]$ such that
there exists a homomorphism of dialgebras
$\varphi: A\to P^{(0)}$ and
$\varphi(a)*\varphi(b) = \sum\limits_i h_i\otimes 1\otimes_H c_i$,
$\deg h_i\le 1$
for all $a,b\in A$.
Then there exists a homomorphism of pseudo-algebras
$C_{\Var}(\varphi): C_{\Var}(A) \to P$
such that $C_{\Var}(\varphi)|_A = \varphi$.
\end{prop}

\begin{proof}
Let us first consider  the case when $\Sigma =\emptyset$,
i.e., $A$ is an arbitrary 0-dialgebra.

Suppose $P$ is a pseudo-algebra over $H=\Bbbk[T]$ that satisfies
the conditions of the proposition.
Since $\Bbbk $ may be of positive characteristic, we need a slightly
different definition of $n$-products:
\[
  x\oo{n} y = z_n,\quad
x,y\in P,
\]
where
\[
  x*y= \sum\limits_{s\ge 0}T^s\otimes 1\otimes _H z_s .
\]
Then
$x\oo{n} y=0$ for $n\gg 0$, and
$H^{\otimes 2}$-linearity implies
\begin{equation}\label{eq:newC}
 Tx\oo{n} y = x\oo{n-1} y, \quad x\oo{n} Ty = T(x\oo{n}y) -x\oo{n-1} y,
\end{equation}
assuming $x\oo{-1} y = 0$.

For all $a,b\in A$ we have
\begin{equation}\label{eq:dhom}
  \varphi(a\rbar b) = \varphi(a)\oo{0} \varphi(b),
\quad
  \varphi(a\lbar b)=  \varphi(a)\oo{0} \varphi(b) +
  T(\varphi(a)\oo{1} \varphi(b))
\end{equation}
since $\varphi : A\to P^{(0)}$ is a homomorphism of dialgebras.

\begin{lem}\label{lem:duo}
For all $a,b,c,d\in A$ we have
$(\varphi(a)\oo{1} \varphi(b))\oo{n} (\varphi(c)\oo{1} \varphi(d))=0$
for $n\ge 0$, and
$(\varphi(a)\oo{1} \varphi(b))\oo{n} \varphi(c) = 0$,
$\varphi(a)\oo{n}(\varphi(b)\oo{1} \varphi(c))=0$
for $n\ge 1$.
\end{lem}

\begin{proof}
Denote $\varphi(a)$, $\varphi(b)$, $\varphi(c)$,
$\varphi(d)$ by $x$, $y$, $z$, $u$, respectively.

There exists $N\ge 0$ such that $x\oo{n} (y\oo{1} z)=0$
for all $n\ge N$. Suppose that the minimal of such $N$ greater than~2.
Then  by \eqref{eq:dhom}
\[
 x\oo{N} T(y\oo{1} z) = x\oo{N}\varphi(b\rbar c)-x\oo{N}\varphi(b\lbar c)=0.
\]
On the other hand, by \eqref{eq:newC}
\[
 x\oo{N} T(y\oo{1} z) = T(x\oo{N} (y\oo{1} z)) - x\oo{N-1} (y\oo{1} z),
\]
so
$x\oo{N-1} (y\oo{1} z)=0$ in contradiction with the minimality of $N$.
Other statements of the Lemma can be proved in a completely similar way.
\end{proof}

Define a map $\psi: C(A) \to P$ by the following rule:
\[
  \psi(h\otimes a) = h(T)\varphi(a),
\quad
  \psi(a\otimes b ) = -\varphi(a)\oo{1} \varphi(b).
\]
The second equation actually means that we define
$\psi: A\otimes A \to P$, then by Lemma~\ref{lem:duo}
 we have $U=\Span\{\langle a,b\rangle\otimes \langle c, d\rangle\mid
a,b,c,d \in A\}\subseteq \Ker\psi$.

It remains to show that $\psi $ is actually a homomorphism
of pseudo-algebras. Relation \eqref{eq:dhom} implies that
$\psi $ is $H$-linear. To check that $\psi $
preserves $*$, one should proceed with four cases
following~\eqref{eq:*1}--\eqref{eq:*4}.
Let us consider the second case for example:
$a*(b\otimes c)= 1\otimes 1\otimes _H a\otimes \langle b,c\rangle$,
$\psi(a\otimes \langle b,c\rangle)
=- \varphi(a)\oo{1} \varphi(\langle b,c\rangle)
=\varphi(a) \oo{1} T(\varphi(b)\oo{1} \varphi(c))
= T(\varphi(a) \oo{1} (\varphi(b)\oo{1} \varphi(c)))
 - \varphi  (a) \oo{0} (\varphi(b)\oo{1} \varphi(c))
= - \varphi  (a) \oo{0} (\varphi(b)\oo{1} \varphi(c))$
by Lemma~\ref{lem:duo}.
On the other hand,
$\psi(a)*\psi(b\otimes c ) =
  -1\otimes 1\otimes_H \varphi(a)\oo{0} (\varphi(b)\oo{1} \varphi(c))$.

Now, let $\Sigma $ be a nonempty family of polylinear identities defining
a homogeneous variety $\Var $,
and let $A$ be a $\Var$-dialgebra. As we have proved,
for any pseudo-algebra $P$ satisfying the conditions of the theorem
there exists a homomorphism of pseudo-algebras
$C(\varphi): C(A) \to P$ such that $C(\varphi)|_A = \varphi$.
Recall that $C_{\Var} (A) = C/I$, where $I$ is spanned by
all coefficients of
$C(A)_n(t^*)(u_1,\dots, u_n)\in H^{\otimes n}\otimes_H C(A)$,
$t\in \Sigma $, $n=\deg t$,
$u_i\in C(A)$, $i=1,\dots, n$.
Therefore, $I\subseteq \Ker C(\varphi)$ and
there exists a homomorphism $C_{\Var }(\varphi)$ with the required
properties.
\end{proof}

\section{Conformal representations of Leibniz algebras}

Recall the notions of a conformal endomorphism
\cite{K1,K2} and a representation
of an (associative or Lie) pseudo-algebra \cite{BDK}.
We will generally follow \cite{Ko4} in the exposition.

Consider a left module $M$ over $H=\Bbbk[T]$.
By $\End M$ we denote the associative algebra of $\Bbbk$-linear
maps from $M$ to itself.
A map $a: \Bbbk \to \End M$ is said to be {\em locally regular\/}
if for any $u\in M$ the map $z\mapsto a(z)u$ is a regular
$M$-valued map on the affine line~$\Bbbk$.

Note that $\End M$ can be considered as a left $H$-module
 with respect to
\[
(f \varphi)(u) = f_{(2)}\varphi(f_{(-1)}u), \quad
   f\in H,\ u\in M,\ \varphi\in \End M.
\]
A map
$a: \Bbbk \to \End M$ is {\em translation-invariant\/}
if
\[
fa(z) = f(z)a(z), \quad f\in A, \ z\in V.
\]

\begin{defn}
A {\em conformal endomorphism\/} of an $H$-module
$M$ is a locally regular and
translation invariant map $a : \Bbbk \to \End M$.
\end{defn}

The set $\Cend M$ of all conformal endomorphisms of $M$
is a linear space
which is also an $H$-module with respect to the action defined by
\[
(fa)(z)= f(z)a(z), \quad f\in H,\ a\in \Cend M,\ z\in \Bbbk.
\]

It is easy to see that an $H$-pseudo-algebra structure on a
left $H$-module $C$ can be equivalently
defined via an $H$-linear map
$\mu : C\to \Cend C$.
Indeed, the correspondence between the map $\mu $
and pseudo-product $*$ is completely described as follows: if
\[
\mu(a)(x)b= \sum\limits_{s=0}^n x^s c_s, \quad b,c_s\in C,
\]
then
\[
a*b = \sum\limits_{s=0}^n (-T)^s \otimes 1 \otimes _H c_s.
\]

For the $H$-module  $\Cend M$ and for every $a\in \Cend M$
one may define a map
$(a_{x}\cdot): \Bbbk \to \End \Cend M$ as follows:
\[
 (a_{z} b)(w) = a(z)b(w-z), \quad b\in \Cend M,\ z,w\in \Bbbk.
\]
This map may not be locally regular, so $\Cend M$
is not necessarily a pseudo-algebra (see \cite{BDK,K1,K2}).
However, some of its submodules are actually pseudo-algebras.
For example, if $M$ is a free $H$-module
generated by a space $M_0$
(i.e., $M=H\otimes M_0$)
then
the current submodule
$\Curr \End M_0 = H\otimes \End M_0 \subset \Cend M$,
where
$(f\otimes \varphi)(z)= f(-z)(\idd_H\otimes \varphi)$,
$f\in H$, $\varphi \in \End M_0$, $z\in \Bbbk $,
is always an associative pseudo-algebra.

Every associative pseudo-algebra $C$ can be considered as a Lie
pseudo-algebra with respect to new pseudo-product
\[
 [a*b] = a*b - ((12)\otimes_H \idd_C)(b*a).
\]
The Lie pseudo-algebra obtained is denoted by $C^{(-)}$.
In particular, if the entire $\Cend M$ is a pseudo-algebra
(that is the case if $M$ is a finitely generated $H$-module)
then $\Cend M^{(-)}$ is denoted by $\gc M$.
Even if $\Cend M$ is not a pseudo-algebra, we will use the following
convention: $L\subset \gc M$ means that $L$ is a subalgebra of
some $C^{(-)}$, where $C\subseteq \Cend M$ is a pseudo-algebra.

\begin{defn}\label{defn:conf.repr}
Let $\mathfrak g$ be a Leibniz algebra, i.e., a Lie dialgebra.
A homomorphism of dialgebras $\rho: \mathfrak g \to L^{(0)}$,
where $L\subseteq \gc M$, is called a {\em conformal
representation\/} of $\mathfrak g$ on~$M$.

A conformal representation is said to be finite if $M$
is a finitely generated $H$-module.
\end{defn}

\begin{prop}\label{prop:Ado}
Every Leibniz algebra $\mathfrak g$
 has a faithful conformal representation.
If $\dim\mathfrak g<\infty $ then there exists a finite
faithful conformal representation.
\end{prop}

\begin{proof}
Let us use $[ab]$ for $a\lbar b$, $a,b\in \mathfrak g$.
Denote by
$\mathfrak l$
the (ordinary) Lie algebra
$\mathfrak g/\Span\{[xx] \mid x\in \mathfrak g\}$.
If $x\in \mathfrak g$ then $\bar x$ is the image of $x$
in~$\mathfrak l$.

Let $V$ stands for an arbitrary nonzero (left) $\mathfrak l$-module.
Set $M_0 = V\oplus (\mathfrak g\otimes V)$,
$M = H\otimes M_0$.
Consider the map
\[
\rho : \mathfrak g \to \Cend M
\]
defined as follows:
\begin{equation}\label{eq:rho}
\begin{gathered}
\rho(x)*u = 1\otimes 1\otimes_H \bar xu - T\otimes 1\otimes _H (x\otimes u), \\
\rho(x)*(a\otimes u) =
  1\otimes 1\otimes_H (a\otimes \bar x u - [ax]\otimes u),
\end{gathered}
\end{equation}
$a,x\in \mathfrak g$, $u\in V$.
It is easy to see that
$\rho (x) \in \Curr \End M_0$. Indeed,
\[
  \rho(x) = \rho_0(x) - T\rho_1(x),
\]
where $\rho_i(x) \in \End M_0$,
$\rho_0(x) = \bar x \oplus (\idd\otimes \bar x - [\cdot x]\otimes \idd)$,
$\rho_1(x)$ maps
$\mathfrak g\otimes V$ to zero and
$u\in V$ to $x\otimes u$.
Therefore, $\rho $ is injective, and $\rho(\mathfrak g)$
generates  an associative pseudo-algebra $C$ that lies in
$\Curr \End M_0$.
Moreover,
\begin{multline}
\rho(a)*\rho(b) = (\rho_0(a)-T\rho_1(a))*(\rho_0(b)-T\rho_1(b)) \\
=
1\otimes 1\otimes_H \rho_0(a)\rho_0(b)
- T\otimes 1\otimes _H \rho_1(a)\rho_0(b)
-1\otimes T\otimes _H \rho_0(a)\rho_1(b)
+ T\otimes T\otimes_H \rho_1(a)\rho_1(b) \\
=
1\otimes 1\otimes _H
  (\rho_0(a)\rho_0(b)- T\rho_0(a)\rho_1(b))
-T\otimes 1\otimes_H (\rho_1(a)\rho_0(b)- \rho_0(a)\rho_1(b)),
\label{eq:rho*}
\end{multline}
since $\rho_1(a)\rho_1(b)=0$.

This is straightforward to check that
\[
[\rho_0(a),\rho_0(b)] = \rho_0([ab]),
\quad
[\rho_1(a), \rho_0(b)]=\rho_1([ab]).
\]
Then \eqref{eq:C0} and \eqref{eq:rho*} imply
\[
 C^{(0)}_2 (x_1\lbar x_2 - x_2\rbar x_1)(\rho(a), \rho(b))
=\rho([ab]), \quad a,b\in \mathfrak g,
\]
that means $\rho $ is a conformal representation of~$\mathfrak g$.

If $\mathfrak g$ is finite-dimensional then
so can be choosen $V$, and hence $M$ is a
finitely generated $H$-module.
\end{proof}

\begin{cor}[\cite{L2}]
Every Leibniz algebra can be embedded into
an associative dialgebra.
\end{cor}

\begin{proof}
Consider a Leibniz algebra $\mathfrak g$ and its
faithful conformal representation $\rho : \mathfrak g \to L\subset \gc M$
as in Proposition~\ref{prop:Ado}.
We have seen that $L$ is actually a Lie subalgebra of $(\Curr \End M_0)^{(-)}$,
so $\mathfrak g \subset L \subset A=(\Curr \End M_0)^{(0)}$.
The latter is an associative dialgebra by
Theorem~\ref{thmVar0}.
\end{proof}

\begin{rem}
If $V=U(\mathfrak l)$ in Proposition~\ref{prop:Ado},
then the associative dialgebra generated in
$(\Cend M)^{(0)}$ by the image of $\mathfrak g$
is isomorphic
to the universal enveloping dialgebra $Ud(\mathfrak g)$
constructed in \cite{L2}. However, $C_{Lie}(\mathfrak g)$
is not embedded into $\Curr \End M_0$ for $M_0 = U(\mathfrak l)\oplus
(\mathfrak g\otimes U(\mathfrak l))$.
\end{rem}


\begin{thebibliography}{99}

\bibitem{L0}
J.-L.~Loday,
Une version non commutative des alg`ebres de Lie:
les alg`ebres de Leibniz,
Enseign. Math., 39  (1993)  269–-293.

\bibitem{L1}
J.-L. Loday, T.~Pirashvili,
Universal envelopping algebras of Leibniz algebras and homology,
Math. Ann. {\bf 296} (1993) 139--158.

\bibitem{L2}
J.-L.~Loday,
Dialgebras,
in: Dialgebras and Related Operads,
Lecture Notes in Mathematics, vol.~1763, pp.~7-66.
Springer Verl., Berlin, 2001.

\bibitem{DL}
D.~Liu,
Steinberg--Leibniz algebras and superalgebras,
J. Algebra  283 (1) (2005) 199--221.

\bibitem{Faulk}
 J.R.~Faulkner,
 Barbilian planes,
 Geom. Dedicata 30 (1989) 125--181.

\bibitem{K1}
 V.G.~Kac,
 Vertex Algebras for Beginners,
 University Lecture Series, vol.~10, AMS, Providence, RI, 1996.

\bibitem{Ro1}
M.~Roitman,
On free conformal and vertex algebras,
J. Algebra 217 (2) (1999) 496--527.

\bibitem{Ko3}
P.S.~Kolesnikov,
Identities of conformal algebras and pseudoalgebras,
Comm. Algebra 34 (6) (2006) 1965--1979.

\bibitem{BDK}
 B.~Bakalov, A.~D'Andrea,  V.G.~Kac,
 Theory of finite pseudoalgebras,
 Adv. Math. 162 (1) (2001) 1--140.

\bibitem{Ko4}
P.S.~Kolesnikov,
Associative algebras related to conformal algebras,
Applied Categorical Structures, to appear.

\bibitem{K2}
 V.G.~Kac,
 Formal distribution algebras and conformal algebras,
 XII-th International Congress in Mathematical Physics (ICMP'97)
 (Brisbane), International Press, Cambridge, MA, 1999, pp.~80--97.


\end{thebibliography}
\end{document}